\documentclass[10pt,letterpaper]{article}

\usepackage[T1]{fontenc}
\usepackage[utf8]{inputenc}
\usepackage[sc]{mathpazo}
\usepackage[scaled=0.88]{helvet}
\usepackage{microtype}
\usepackage[letterpaper,left=0.68in,right=0.68in,top=0.60in,bottom=0.56in,headheight=18pt]{geometry}

\usepackage{amsmath,amssymb,amsthm,mathtools}
\usepackage{booktabs,tabularx,array}
\usepackage{xcolor,colortbl}
\usepackage[most]{tcolorbox}
\usepackage{fancyhdr}
\usepackage{titlesec}
\usepackage{enumitem}
\usepackage[hypcap=false]{caption}
\usepackage{tikz}
\usepackage{pgfplots}
\pgfplotsset{compat=1.18}
\usetikzlibrary{arrows.meta,calc}

\usepackage[colorlinks=true,linkcolor=primary,citecolor=slate,urlcolor=slate]{hyperref}

\newcommand{\ml}{\mathcal{L}}
\newcommand{\md}{\mathcal{D}}
\newcommand{\dc}{d_{\mathrm{critical}}}
\newcommand{\E}{\mathbf{E}}
\newcommand{\Wm}{W_{-1}}

\definecolor{primary}{RGB}{22, 56, 92}
\definecolor{slate}{RGB}{42, 102, 158}
\definecolor{okgreen}{RGB}{18, 122, 74}
\definecolor{warnamber}{RGB}{186, 108, 10}
\definecolor{crimson}{RGB}{176, 42, 46}
\definecolor{muted}{RGB}{90, 100, 112}
\definecolor{rulecol}{RGB}{206, 214, 224}
\definecolor{tintbg}{RGB}{244, 247, 251}
\definecolor{greentint}{RGB}{238, 248, 243}
\definecolor{ambertint}{RGB}{254, 247, 236}

\setlist[itemize]{leftmargin=1.4em,itemsep=1.5pt,topsep=1.5pt}
\setlist[enumerate]{leftmargin=1.6em,itemsep=1.5pt,topsep=1.5pt}

\titleformat{\section}{\large\sffamily\bfseries\color{primary}}{\thesection}{0.65em}{}[\vspace{2pt}{\color{rulecol}\titlerule[0.75pt]}]
\titleformat{\subsection}{\normalsize\sffamily\bfseries\color{slate}}{\thesubsection}{0.55em}{}
\titlespacing*{\section}{0pt}{9pt}{4pt}
\titlespacing*{\subsection}{0pt}{6pt}{2pt}

\renewcommand{\headrulewidth}{0.4pt}
\renewcommand{\headrule}{\hbox to\headwidth{\color{rulecol}\leaders\hrule height \headrulewidth\hfill}}

\theoremstyle{plain}
\newtheorem{theorem}{Theorem}
\newtheorem{lemma}{Lemma}
\newtheorem{corollary}{Corollary}
\theoremstyle{definition}
\newtheorem{remark}{Remark}

\newtcolorbox{thmBox}[1][]{
  enhanced, colback=tintbg, colframe=primary,
  boxrule=0.6pt, leftrule=3pt, arc=2pt,
  left=11pt, right=11pt, top=7pt, bottom=7pt,
  fonttitle=\small\sffamily\bfseries\color{primary},
  colbacktitle=tintbg, title=#1, attach title to upper={\par\vspace{3pt}}
}
\newtcolorbox{noteBox}[1][]{
  enhanced, colback=ambertint, colframe=warnamber,
  boxrule=0.5pt, leftrule=3pt, arc=2pt,
  left=9pt, right=9pt, top=5pt, bottom=5pt,
  fontupper=\small,
  fonttitle=\footnotesize\sffamily\bfseries\color{warnamber},
  colbacktitle=ambertint, title=#1, attach title to upper={\par\vspace{1.5pt}}
}

\begin{document}
\thispagestyle{empty}

\begin{tcolorbox}[enhanced,colback=primary,colframe=primary,arc=3pt,left=16pt,right=16pt,top=12pt,bottom=12pt]
  {\footnotesize\sffamily\color{white!80!slate}\MakeUppercase{Technical Note \quad\textbar\quad Companion Proof to Kadloor \& Santhi (2010)}}\\[5pt]
  {\LARGE\bfseries\color{white} Closed-Form Formula for $\dc(a_0)$ in Large\\Fully Connected Power Grids}\\[5pt]
  {\small\sffamily\color{white!86!slate} An exact algebraic solution of the recursive system in \texttt{arXiv:1011.4098v1} (Section~4, \texttt{0802.tex})}
\end{tcolorbox}

\vspace{2pt}
\textbf{\sffamily\small\color{primary}Scope.}\enspace
In Section~4 of S.\ Kadloor and N.\ Santhi (2010), \textit{``Understanding Cascading Failures in Power Grids''} (source \texttt{0802.tex} in \texttt{arXiv:1011.4098v1}), the evolution of a cascading failure on a large fully connected network ($N \to \infty$) is given by a four-variable recursive system in the sequences $(a_n, \mu_{n-1}, \md_n, p_n)$. This self-contained note uses the paper's \emph{exact} notation, symbols, and indexing throughout to prove that the coupled recursion telescopes into a single one-dimensional map $a_{n+1} = g(a_n)$, yielding a closed-form expression for the critical mean disturbance $\dc(a_0)$ in terms of the lower real branch $\Wm$ of the Lambert~$W$ function.

\section{System Model, Notation, and Recursive Equations}
\label{sec:setup}

We consider the homogeneous large fully connected network $G(V,E)$ of $N = |V|$ generator nodes from Sections~2 and~4 of~\cite{ks2010}: every node $v \in V$ has maximum load capacity $L_{\max} = 1$, constant initial load $a_0 \in [0, 1)$, and receives an independent additive disturbance $D_v(0) \sim \mathrm{Exponential}(d_m)$ of mean $d_m > 0$ at time~$0$. Table~\ref{tab:notation} aligns every symbol directly with \texttt{0802.tex}.

\begin{table}[!htbp]
\centering
\caption{Symbol glossary matching \texttt{0802.tex} (\texttt{arXiv:1011.4098v1}, Sections~2 and~4).}
\label{tab:notation}
\small
\begin{tabularx}{\linewidth}{@{}l>{\raggedright\arraybackslash}lX@{}}
\toprule
Symbol & Paper location & Meaning \\
\midrule
$L_{\max}$ & Sec.~2 & Maximum capacity of each generator node ($L_{\max} = 1$ unit; general $L_{\max} > 0$ in Corollary~\ref{cor:scale}). \\
$a_0$ & Sec.~2 & Initial steady-state load at every node, $a_0 \in [0, 1)$ (paper focuses on $a_0 \in (0, 1)$). \\
$\ml_0$ & Sec.~4 & Random variable of the initial load: $\ml_0 \sim \delta(l - a_0)$ (constant at $a_0$). \\
$\md_0$ & Sec.~4 & Initial disturbance random variable: $\md_0 \sim \mathrm{Exponential}(d_m)$, density $\tfrac{1}{d_m} e^{-x/d_m}$ ($x \ge 0$). \\
$d_m$ & Sec.~2 & Mean value of the initial exponential disturbance $\md_0$ ($d_m > 0$; \cite{ks2010} assumes $d_m < 1$). \\
$\ml_n$ & Sec.~4 & Random variable with the distribution of load at surviving nodes at stage $n$. \\
$\md_n$ & Sec.~4 & Deterministic load redistributed equally to each surviving node at stage $n \ge 1$. \\
$p_n$ & Sec.~4 & Conditional probability that a node surviving into stage $n$ fails at the end of stage $n$. \\
$\mu_n$ & Sec.~4 & Mean load $\E(L_{k_i}(n))$ at nodes $k_i$ that fail at the end of stage $n$ ($n \ge 0$). \\
$a_n$ & Sec.~4 & Lower endpoint (minimum shifted load) of the surviving density $\ml_n$ ($a_1 = a_0$, $a_n = a_{n-1} + \md_{n-1}$). \\
$\dc$ & Sec.~2 & Critical threshold for $d_m$: cascade stops if $d_m \le \dc$, complete outage if $d_m > \dc$. \\
\bottomrule
\end{tabularx}
\end{table}

\paragraph{Recursive equations (\texttt{0802.tex}, Section~4, lines~211--222).}
In the limit $N \to \infty$ (via Lemma~4.1 of~\cite{ks2010}), the cascade is governed verbatim by:
\begin{itemize}
  \item \textbf{Initialization ($n = 0, 1$):} with $\mu_0 = \E(L_{k_i}(0)) = 1 + d_m$,
  \begin{equation}
  \label{eq:init}
    p_0 = e^{-\frac{1 - a_0}{d_m}}, \qquad
    \md_1 = \frac{p_0}{1 - p_0}\,(1 + d_m), \qquad
    a_1 = a_0, \qquad
    p_1 = \frac{e^{-\frac{1 - a_1}{d_m}}}{1 - e^{-\frac{1 - a_1}{d_m}}}\left(e^{\frac{\md_1}{d_m}} - 1\right).
  \end{equation}
  \item \textbf{Iteration ($n = 2, 3, \dots$):} If $\md_{n-1} \ge 1 - a_{n-1}$ (equivalently $a_{n-1} + \md_{n-1} \ge 1$), stop (all remaining nodes exceed capacity $L_{\max}=1$ and the entire grid fails). Otherwise, compute:
  \begin{subequations}
  \label{eq:rec}
  \begin{align}
    a_n &= a_{n-1} + \md_{n-1}, \label{eq:rec-a} \\
    \mu_{n-1} &= 1 + d_m - \frac{\md_{n-1}}{e^{\md_{n-1}/d_m} - 1}, \label{eq:rec-mu} \\
    \md_n &= \frac{p_{n-1}}{1 - p_{n-1}}\,\mu_{n-1}, \label{eq:rec-D} \\
    p_n &= \frac{e^{-\frac{1 - a_n}{d_m}}}{1 - e^{-\frac{1 - a_n}{d_m}}}\left(e^{\frac{\md_n}{d_m}} - 1\right). \label{eq:rec-p}
  \end{align}
  \end{subequations}
\end{itemize}
\textit{Note on notation:} In \eqref{eq:init}--\eqref{eq:rec}, each term $\md_n$ ($n \ge 1$) is strictly positive while $a_n < 1$, and by \eqref{eq:rec-a} the shifted minimum load satisfies $a_n = a_0 + \sum_{k=1}^{n-1} \md_k$. Thus the cumulative load added to every surviving node before stage~$n$ is simply $a_n - a_0$.

\begin{thmBox}[Main Theorem (Exact Analytical Formula for $\dc$)]
\begin{theorem}
\label{thm:main}
Fix $a_0 \in [0, 1)$ and consider the recursive system \eqref{eq:init}--\eqref{eq:rec} with mean initial disturbance $d_m > 0$. There exists a unique constant $\dc = \dc(a_0) \in (0, 1 - a_0)$ defined implicitly by
\begin{equation}
\label{eq:implicit}
  1 - a_0 \;=\; \dc \left(1 + \ln\!\left(1 + \frac{1}{\dc}\right)\right),
\end{equation}
and explicitly in closed form via the lower real branch $\Wm : (-e^{-1}, 0) \to (-\infty, -1)$ of the Lambert~$W$ function:
\begin{equation}
\label{eq:lambert}
  \dc(a_0) \;=\; \frac{1 - a_0}{-\Wm\!\bigl(-(1 - a_0)\,e^{-(2 - a_0)}\bigr) - (1 - a_0)}.
\end{equation}
The sequence $\{a_n\}$ exhibits a sharp threshold dichotomy at $\dc(a_0)$:
\begin{enumerate}[label=(\roman*)]
  \item \textbf{Subcritical \& critical survival ($d_m \le \dc(a_0)$):} The cascade never hits the stopping guard ($\md_{n-1} < 1 - a_{n-1}$ for all $n \ge 2$). The sequence $\{a_n\}$ increases strictly to a limit $a^* = \lim_{n \to \infty} a_n \in (a_0, 1)$, with $\lim_{n\to\infty}\md_n = 0$ and $\lim_{n\to\infty} p_n = 0$. At criticality ($d_m = \dc$), the limiting minimum load is exactly $a^* = a_0 + \dc$.
  \item \textbf{Supercritical blackout ($d_m > \dc(a_0)$):} In finitely many steps $n \ge 2$, the sequence triggers the stopping guard $\md_{n-1} \ge 1 - a_{n-1}$ ($a_{n-1} + \md_{n-1} \ge 1$ and $p_{n-1} \ge 1$), so all surviving nodes overload and the network suffers a complete outage.
\end{enumerate}
\end{theorem}
\end{thmBox}

\section{Proof of Theorem~\ref{thm:main}}
\label{sec:proof}

The proof proceeds directly from \eqref{eq:init}--\eqref{eq:rec} in three short steps: (1) an exact algebraic telescoping identity that eliminates $(\mu_{n-1}, \md_n, p_n)$ without any heuristic approximation; (2) monotonicity and concavity of the resulting one-dimensional iteration $a_{n+1} = g(a_n)$; and (3) uniqueness and inversion of \eqref{eq:implicit} via the Lambert~$W$ function.

For compactness, for any $a \in [a_0, 1)$ and $d_m > 0$ define
\[
  r(a) \;:=\; e^{-\frac{1 - a}{d_m}} \in (0, 1), \qquad
  h(a) \;:=\; (a - a_0) - (1 + d_m)\,r(a), \qquad
  g(a) \;:=\; \frac{a_0 + (1 + d_m - a)\,r(a)}{1 - r(a)}.
\]
Note that $r(a_0) = p_0$, and algebraically
\begin{equation}
\label{eq:g-shift}
  g(a) - a \;=\; \frac{(1 + d_m)\,r(a) - (a - a_0)}{1 - r(a)} \;=\; \frac{-h(a)}{1 - r(a)}.
\end{equation}

\begin{lemma}[Exact telescoping identity and reduction to a 1D map]
\label{lem:telescope}
Along the trajectory of \eqref{eq:init}--\eqref{eq:rec}, as long as $a_n < 1$ (i.e.\ $\md_{k-1} < 1 - a_{k-1}$ for all $2 \le k \le n$), we have $a_n \in [a_0, 1)$, $\md_n > 0$, and $p_{n-1} \in (0, 1)$, and the iterates satisfy
\begin{equation}
\label{eq:inv}
  (a_n - a_0) + \md_n \;=\; \bigl(1 + d_m + \md_n\bigr)\,e^{-\frac{1 - a_n}{d_m}} \qquad (n \ge 1).
\end{equation}
Equivalently, $\md_n = \bigl[(1+d_m)\,r(a_n) - (a_n - a_0)\bigr] / \bigl[1 - r(a_n)\bigr]$, so if the update continues to stage $n+1$ then
\begin{equation}
\label{eq:1dmap}
  a_{n+1} \;=\; a_n + \md_n \;=\; g(a_n) \qquad (n \ge 1).
\end{equation}
Moreover, the stagewise survival probabilities telescope: for every $k \ge 1$ with $a_{k+1} < 1$,
\begin{equation}
\label{eq:surv-tel}
  1 - p_k \;=\; \frac{1 - r(a_{k+1})}{1 - r(a_k)}, \qquad\text{and hence}\qquad
  \prod_{k=0}^{n-1} (1 - p_k) \;=\; 1 - e^{-\frac{1 - a_n}{d_m}}.
\end{equation}
\end{lemma}

\begin{proof}
We proceed by induction on $n \ge 1$. For the base case $n = 1$, $a_1 = a_0$ gives $a_1 - a_0 = 0$ and $r(a_1) = p_0 \in (0, 1)$. From \eqref{eq:init}, $(1 - r(a_1))\,\md_1 = (1 + d_m)\,r(a_1)$, which rearranges immediately to $\md_1 = (1 + d_m + \md_1)\,r(a_1)$, verifying \eqref{eq:inv} for $n=1$ (with $\md_1 > 0$).

Now let $n \ge 2$ and assume \eqref{eq:inv} holds at index $n-1$ and that $a_n = a_{n-1} + \md_{n-1} < 1$ with $\md_{n-1} > 0$. Because $r(a_{n-1})\,e^{\md_{n-1}/d_m} = e^{-(1 - a_{n-1} - \md_{n-1})/d_m} = e^{-(1 - a_n)/d_m} = r(a_n) < 1$, substituting \eqref{eq:rec-p} at index $n-1$ yields
\begin{equation}
\label{eq:1mp}
  1 - p_{n-1} \;=\; 1 - \frac{r(a_{n-1})\bigl(e^{\md_{n-1}/d_m} - 1\bigr)}{1 - r(a_{n-1})}
  \;=\; \frac{1 - r(a_{n-1})\,e^{\md_{n-1}/d_m}}{1 - r(a_{n-1})}
  \;=\; \frac{1 - r(a_n)}{1 - r(a_{n-1})} \;>\; 0,
\end{equation}
which establishes the first part of \eqref{eq:surv-tel} (and ensures $p_{n-1} \in (0, 1)$ since $\md_{n-1} > 0$). Multiplying \eqref{eq:1mp} over $k = 1, \dots, n-1$ together with $1 - p_0 = 1 - r(a_1)$ proves the cumulative survival product in \eqref{eq:surv-tel}.

From \eqref{eq:rec-p} and \eqref{eq:1mp}, the denominator $1 - r(a_{n-1})$ cancels in the ratio
\[
  \frac{p_{n-1}}{1 - p_{n-1}} \;=\; \frac{r(a_{n-1})\bigl(e^{\md_{n-1}/d_m} - 1\bigr)}{1 - r(a_n)}.
\]
Inserting this and \eqref{eq:rec-mu} into \eqref{eq:rec-D} cancels the factor $(e^{\md_{n-1}/d_m} - 1)$ identically:
\begin{align*}
  \bigl(1 - r(a_n)\bigr)\,\md_n
  &= r(a_{n-1})\bigl(e^{\md_{n-1}/d_m} - 1\bigr)\!\left(1 + d_m - \frac{\md_{n-1}}{e^{\md_{n-1}/d_m} - 1}\right) \\
  &= (1 + d_m)\bigl(r(a_n) - r(a_{n-1})\bigr) - \md_{n-1}\,r(a_{n-1}).
\end{align*}
Note that since $e^x - 1 > x$ for all $x > 0$, setting $x = \md_{n-1}/d_m$ gives $\mu_{n-1} > 1 > 0$, and hence $\md_n > 0$. Moving $-r(a_n)\,\md_n$ to the right-hand side separates the indices cleanly into a telescoping difference:
\[
  \md_n \;=\; \bigl(1 + d_m + \md_n\bigr)\,r(a_n) - \bigl(1 + d_m + \md_{n-1}\bigr)\,r(a_{n-1}).
\]
Adding the induction hypothesis $a_n - a_0 = (a_{n-1} - a_0) + \md_{n-1} = (1 + d_m + \md_{n-1})\,r(a_{n-1})$ cancels the entire $(n-1)$ term, yielding $(a_n - a_0) + \md_n = (1 + d_m + \md_n)\,r(a_n)$, which is \eqref{eq:inv}. Solving \eqref{eq:inv} for $\md_n$ gives $\md_n = -h(a_n)/(1 - r(a_n))$, so $a_n + \md_n = g(a_n)$ by \eqref{eq:g-shift}.
\end{proof}

\begin{remark}[Physical interpretation via total load conservation]
Identity \eqref{eq:inv} directly encodes conservation of load on the fully connected network $G(V,E)$. Prior to stage~$n$, every surviving node carries $a_0 + \md_0 + (a_n - a_0) = a_n + \md_0 < 1$, so the nodes that have already failed are exactly those with initial disturbance $\md_0 \ge 1 - a_n$ (fraction $r(a_n) = e^{-(1-a_n)/d_m}$). Since all redistributed load originates from the initial loads of earlier failed nodes (no double-counting), summing only the initial loads over $\md_0 \ge 1 - a_n$ via the memoryless property of $\mathrm{Exponential}(d_m)$ gives $\E[(a_0 + \md_0)\,\mathbf{1}_{\md_0 \ge 1 - a_n}] = (a_0 + 1 - a_n + d_m)\,r(a_n)$ (where $\mathbf{1}_{\md_0 \ge 1 - a_n}$ is the indicator random variable of the event $\md_0 \ge 1 - a_n$). Distributing this load equally across the surviving fraction $1 - r(a_n)$ yields the next cumulative shift $a_{n+1} - a_0 = (1 - a_n + a_0 + d_m)\,r(a_n) / (1 - r(a_n))$, which rearranges identically to \eqref{eq:inv}.
\end{remark}

\begin{lemma}[Strict monotonicity of $g$ and concavity of $h$]
\label{lem:mono}
For every $a_0 \in [0, 1)$ and $d_m > 0$:
\begin{enumerate}[label=(\alph*)]
  \item $g$ is strictly increasing on $[a_0, 1)$, with derivative
  \[
    g'(a) \;=\; \frac{r(a)}{d_m\,\bigl(1 - r(a)\bigr)^2}\,\Bigl[(1 - a) + a_0 + d_m\,r(a)\Bigr] \;>\; 0.
  \]
  \item $h(a) = (a - a_0) - (1 + d_m)\,r(a)$ is strictly concave on $\mathbb{R}$, with unique global maximizer
  \[
    a_{\mathrm{peak}} \;=\; 1 - d_m \ln\!\left(1 + \frac{1}{d_m}\right) \;<\; 1, \qquad
    h(a_{\mathrm{peak}}) \;=\; (1 - a_0) - d_m\!\left(1 + \ln\!\left(1 + \frac{1}{d_m}\right)\right).
  \]
\end{enumerate}
\end{lemma}

\begin{proof}
(a) Since $r'(a) = r(a)/d_m$, write $g(a) = u(a) / v(a)$ with $u(a) = a_0 + (1 + d_m - a)\,r(a)$ and $v(a) = 1 - r(a)$. Then $u'(a) = -r(a) + (1 + d_m - a)\frac{r(a)}{d_m} = \frac{r(a)}{d_m}(1 - a)$ and $v'(a) = -\frac{r(a)}{d_m}$. Applying the quotient rule:
\begin{align*}
  u'(a)v(a) - u(a)v'(a)
  &= \frac{r(a)}{d_m}\Bigl[(1 - a)\bigl(1 - r(a)\bigr) + a_0 + (1 + d_m - a)\,r(a)\Bigr] \\
  &= \frac{r(a)}{d_m}\Bigl[(1 - a) + a_0 + d_m\,r(a)\Bigr].
\end{align*}
For $a \in [a_0, 1)$, we have $1 - a > 0$, $a_0 \ge 0$, $r(a) > 0$, and $d_m > 0$, so $g'(a) > 0$.

(b) Differentiating $h$ twice gives
\[
  h'(a) = 1 - \frac{1 + d_m}{d_m}\,e^{-\frac{1 - a}{d_m}}, \qquad
  h''(a) = -\frac{1 + d_m}{d_m^2}\,e^{-\frac{1 - a}{d_m}} < 0 \quad (\forall\, a \in \mathbb{R}).
\]
Solving $h'(a_{\mathrm{peak}}) = 0$ yields $e^{-(1 - a_{\mathrm{peak}})/d_m} = \frac{d_m}{1 + d_m}$, i.e., $a_{\mathrm{peak}} = 1 - d_m \ln(1 + 1/d_m) < 1$. At this point $(1 + d_m)\,r(a_{\mathrm{peak}}) = d_m$, so $h(a_{\mathrm{peak}}) = a_{\mathrm{peak}} - a_0 - d_m = (1 - a_0) - d_m\bigl(1 + \ln(1 + 1/d_m)\bigr)$.
\end{proof}

\begin{lemma}[Existence, uniqueness, and Lambert $W_{-1}$ inversion]
\label{lem:lambert}
For $d > 0$, define $\psi(d) := d\bigl(1 + \ln(1 + 1/d)\bigr)$. Then $\psi$ is strictly increasing on $(0, \infty)$ with $\lim_{d\to 0^+} \psi(d) = 0$ and $\psi(1) = 1 + \ln 2 > 1$. Hence for each $a_0 \in [0, 1)$, equation \eqref{eq:implicit} ($\psi(\dc) = 1 - a_0$) has a unique solution $\dc \in (0, 1 - a_0)$, given explicitly by \eqref{eq:lambert}.
\end{lemma}

\begin{proof}
Differentiating $\psi(d)$ yields
\[
  \psi'(d) \;=\; 1 + \ln\!\left(1 + \frac{1}{d}\right) - \frac{1}{1 + d}
  \;=\; \ln\!\left(1 + \frac{1}{d}\right) + \frac{d}{1 + d} \;>\; 0 \quad (\forall\, d > 0).
\]
Since $\psi(d) > d$ for all $d > 0$, $\psi(d) = 1 - a_0$ forces $d < 1 - a_0$. Also $\lim_{d \to 0^+} d\ln(1 + 1/d) = 0$, and $\psi(1) = 1 + \ln 2 \approx 1.693 > 1 \ge 1 - a_0$. By the intermediate value theorem and strict monotonicity, a unique root $\dc \in (0, 1 - a_0)$ exists.

To obtain \eqref{eq:lambert}, put $c := 1 - a_0 \in (0, 1]$ and $u := 1 + 1/\dc > 1$, so $\dc = (u - 1)^{-1}$. Then
\[
  \psi(\dc) = c
  \iff \frac{1 + \ln u}{u - 1} = c
  \iff \ln u - c u = -(1 + c).
\]
Exponentiating both sides and multiplying by $-c$ yields $(-c u)\,e^{-cu} = -c\,e^{-(1+c)}$. Set $w := -cu$. Because $u > 1$, we have $\ln u > 0$, hence $cu = 1 + c + \ln u > 1 + c > 1$, so $w < -1$. Moreover, on $(0, 1]$ the map $c \mapsto c\,e^{-c}$ attains its maximum at $c = 1$, giving
\[
  -c\,e^{-(1+c)} \;=\; -e^{-1}\bigl(c\,e^{-c}\bigr) \;\in\; \bigl[-e^{-2},\, 0\bigr) \;\subset\; \bigl(-e^{-1},\, 0\bigr).
\]
Thus $w$ lies strictly on the lower real branch of the Lambert $W$ function: $w = \Wm\!\bigl(-c\,e^{-(1+c)}\bigr)$. Substituting $1 + 1/\dc = u = -w/c$ yields $\dc = c / (-w - c)$, which is \eqref{eq:lambert}.
\end{proof}

\begin{proof}[Proof of Theorem~\ref{thm:main}]
Set $\dc = \dc(a_0)$ from Lemma~\ref{lem:lambert}. By Lemma~\ref{lem:mono}(b), $h(a_{\mathrm{peak}}) = (1 - a_0) - \psi(d_m)$, so $h(a_{\mathrm{peak}}) \ge 0$ if and only if $d_m \le \dc$ (since $\psi$ is strictly increasing).

\textbf{Case (i): $d_m \le \dc(a_0)$.}\enspace
Here $h(a_{\mathrm{peak}}) \ge 0$, which implies $a_{\mathrm{peak}} = a_0 + d_m + h(a_{\mathrm{peak}}) > a_0$. At the endpoints of $[a_0, a_{\mathrm{peak}}]$, we have $h(a_0) = -(1 + d_m)\,r(a_0) < 0$ and $h(a_{\mathrm{peak}}) \ge 0$. Because $h$ is strictly increasing and continuous on $[a_0, a_{\mathrm{peak}}]$, there exists a \emph{unique} smallest root $a^* \in (a_0, a_{\mathrm{peak}}]$ satisfying $h(a^*) = 0$, i.e., $g(a^*) = a^*$. (Note that $a^* < 1$ because $a_{\mathrm{peak}} < 1$; when $d_m = \dc$, $h(a_{\mathrm{peak}}) = 0$ gives $a^* = a_{\mathrm{peak}} = a_0 + \dc$.)

On $[a_0, a^*)$, $h(a) < 0$, so by \eqref{eq:g-shift} we have $g(a) > a$. We show by induction that $a_n \in [a_0, a^*)$ and $a_{n+1} > a_n$ for every $n \ge 1$. For $n = 1$, $a_1 = a_0 < a^*$. Assuming $a_n \in [a_0, a^*)$, strict monotonicity of $g$ (Lemma~\ref{lem:mono}(a)) gives
\[
  a_n \;<\; g(a_n) \;<\; g(a^*) \;=\; a^* \;<\; 1.
\]
Hence $\md_n = g(a_n) - a_n < 1 - a_n$, so the stopping condition in \eqref{eq:rec} is \emph{never} triggered and $a_{n+1} = g(a_n) \in (a_n, a^*)$. By the monotone convergence theorem, $a_\infty := \lim_{n\to\infty} a_n$ exists and lies in $(a_0, a^*]$. By continuity of $g$ on $[a_0, a^*]$, $a_\infty = g(a_\infty)$, i.e.\ $h(a_\infty) = 0$. Since $a^*$ is the unique root of $h$ in $[a_0, a^*]$, we conclude $a_\infty = a^*$. Consequently, $\md_n = a_{n+1} - a_n \to 0$, and from \eqref{eq:rec-p},
\[
  \lim_{n \to \infty} p_n \;=\; \frac{r(a^*)}{1 - r(a^*)}\bigl(e^0 - 1\bigr) \;=\; 0.
\]

\textbf{Case (ii): $d_m > \dc(a_0)$.}\enspace
Here $h(a_{\mathrm{peak}}) = (1 - a_0) - \psi(d_m) < 0$. By strict concavity of $h$ on $\mathbb{R}$,
\[
  \max_{a \in [a_0, 1]} h(a) \;=\; -\varepsilon \;<\; 0, \qquad\text{where}\qquad
  \varepsilon \;:=\; \begin{cases}
    -h(a_{\mathrm{peak}}) & \text{if } a_{\mathrm{peak}} \ge a_0, \\
    -h(a_0) = (1 + d_m)\,r(a_0) & \text{if } a_{\mathrm{peak}} < a_0.
  \end{cases}
\]
Whenever $a_n \in [a_0, 1)$, Lemma~\ref{lem:telescope} and \eqref{eq:g-shift} imply
\[
  \md_n \;=\; \frac{-h(a_n)}{1 - r(a_n)} \;\ge\; -h(a_n) \;\ge\; \varepsilon \;>\; 0.
\]
Therefore $a_{n+1} = a_n + \md_n \ge a_0 + n\,\varepsilon$ as long as the iteration continues. For some stage $k \le N_{\max} := \lceil (1 - a_0)/\varepsilon \rceil$, we must hit $\md_k \ge 1 - a_k$, triggering the stopping guard in \eqref{eq:rec} at iteration $n = k + 1$ and causing a complete outage in finitely many steps.
\end{proof}

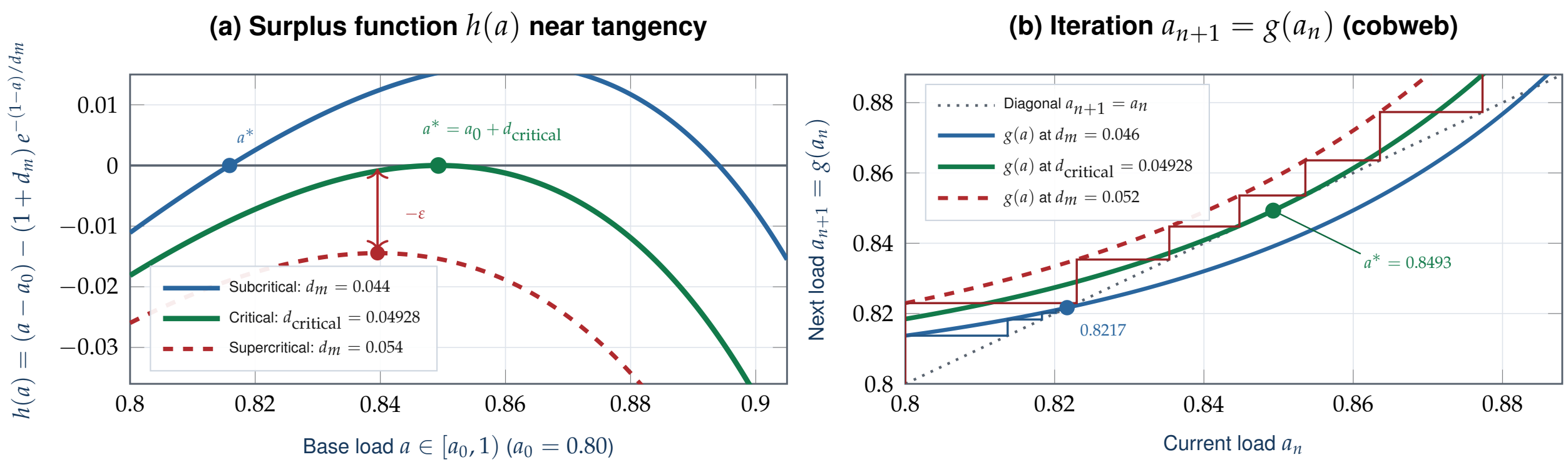
\begin{figure}[!htbp]
\centering
\begin{tikzpicture}
  \begin{axis}[
    name=hplot,
    width=0.47\linewidth, height=4.85cm,
    xmin=0.80, xmax=0.905, ymin=-0.036, ymax=0.015,
    xlabel={Base load $a \in [a_0, 1)$ ($a_0 = 0.80$)},
    ylabel={$h(a) = (a - a_0) - (1 + d_m)\,e^{-(1-a)/d_m}$},
    title={\small\sffamily\bfseries (a) Surplus function $h(a)$ near tangency},
    grid=major, major grid style={draw=rulecol!65,line width=0.3pt},
    axis line style={draw=muted,line width=0.55pt},
    tick label style={font=\scriptsize\sffamily},
    label style={font=\scriptsize\sffamily\color{primary}},
    scaled y ticks=false,
    yticklabel style={/pgf/number format/fixed,/pgf/number format/precision=2},
    ytick={-0.03,-0.02,-0.01,0,0.01},
    legend style={font=\tiny\sffamily, at={(0.03,0.04)}, anchor=south west, fill=white, fill opacity=0.92, draw=rulecol, cells={anchor=west}},
    enlargelimits=false
  ]
    \draw[color=muted, line width=0.7pt] (axis cs:0.80,0) -- (axis cs:0.905,0);
    \addplot[thick, line width=1.4pt, color=slate] table[x=a, y=h_sub] {data/surplus_h.dat};
    \addlegendentry{Subcritical: $d_m = 0.044$}
    \addplot[thick, line width=1.6pt, color=okgreen] table[x=a, y=h_crit] {data/surplus_h.dat};
    \addlegendentry{Critical: $\dc = 0.04928$}
    \addplot[dashed, thick, line width=1.3pt, color=crimson] table[x=a, y=h_sup] {data/surplus_h.dat};
    \addlegendentry{Supercritical: $d_m = 0.054$}

    \addplot[only marks, mark=*, mark size=2.1pt, color=slate] coordinates {(0.81591,0)};
    \addplot[only marks, mark=*, mark size=2.3pt, color=okgreen] coordinates {(0.84928,0)};
    \addplot[only marks, mark=*, mark size=2.0pt, color=crimson] coordinates {(0.83955,-0.01445)};
    \draw[color=crimson, <->, line width=0.7pt] (axis cs:0.83955,-0.0005) -- (axis cs:0.83955,-0.0140);
    \node[font=\tiny\sffamily, color=crimson, anchor=west] at (axis cs:0.842,-0.008) {$-\varepsilon$};
    \node[font=\tiny\sffamily, color=okgreen, anchor=south] at (axis cs:0.858,0.0022) {$a^* = a_0 + \dc$};
    \node[font=\tiny\sffamily, color=slate, anchor=south east] at (axis cs:0.822,0.0018) {$a^*$};
  \end{axis}

  \begin{axis}[
    at={($(hplot.east)+(1.25cm,0)$)}, anchor=west,
    width=0.47\linewidth, height=4.85cm,
    xmin=0.80, xmax=0.888, ymin=0.80, ymax=0.888,
    xlabel={Current load $a_n$},
    ylabel={Next load $a_{n+1} = g(a_n)$},
    title={\small\sffamily\bfseries (b) Iteration $a_{n+1} = g(a_n)$ (cobweb)},
    grid=major, major grid style={draw=rulecol!65,line width=0.3pt},
    axis line style={draw=muted,line width=0.55pt},
    tick label style={font=\scriptsize\sffamily},
    label style={font=\scriptsize\sffamily\color{primary}},
    legend style={font=\tiny\sffamily, at={(0.03,0.97)}, anchor=north west, fill=white, fill opacity=0.92, draw=rulecol, cells={anchor=west}},
    enlargelimits=false
  ]
    \addplot[dotted, thick, color=muted, domain=0.80:0.888] {x};
    \addlegendentry{Diagonal $a_{n+1} = a_n$}
    \addplot[thick, line width=1.2pt, color=slate] table[x=a, y=g_sub] {data/gmap.dat};
    \addlegendentry{$g(a)$ at $d_m = 0.046$}
    \addplot[thick, line width=1.4pt, color=okgreen] table[x=a, y=g_crit] {data/gmap.dat};
    \addlegendentry{$g(a)$ at $\dc = 0.04928$}
    \addplot[dashed, thick, line width=1.2pt, color=crimson] table[x=a, y=g_sup] {data/gmap.dat};
    \addlegendentry{$g(a)$ at $d_m = 0.052$}

    \addplot[color=slate!85!black, line width=0.65pt] table[x=x, y=y] {data/cobweb_sub.dat};
    \addplot[color=crimson!85!black, line width=0.65pt] table[x=x, y=y] {data/cobweb_sup.dat};

    \addplot[only marks, mark=*, mark size=2.1pt, color=slate] coordinates {(0.82167,0.82167)};
    \addplot[only marks, mark=*, mark size=2.2pt, color=okgreen] coordinates {(0.84928,0.84928)};
    \node[font=\tiny\sffamily, color=slate, fill=white, fill opacity=0.88, text opacity=1, inner sep=0.8pt, anchor=north west] at (axis cs:0.823,0.818) {$0.8217$};
    \draw[color=okgreen!85!black, line width=0.45pt] (axis cs:0.84928,0.84928) -- (axis cs:0.861,0.837);
    \node[font=\tiny\sffamily, color=okgreen, fill=white, fill opacity=0.92, text opacity=1, inner sep=0.8pt, anchor=north west] at (axis cs:0.861,0.838) {$a^* = 0.8493$};
  \end{axis}
\end{tikzpicture}
\caption{Geometric picture of the proof at $a_0 = 0.80$. (a) The fixed points of $g$ are zeros of the strictly concave surplus function $h(a)$. Below $\dc$ ($d_m = 0.044$), $h$ crosses zero at $a^* = 0.8159$; at $\dc = 0.04928$, $h$ kisses zero tangentially at $a^* = a_0 + \dc = 0.8493$; above $\dc$ ($d_m = 0.054$), $h(a) \le -\varepsilon < 0$ everywhere. (b) Cobweb diagrams of $a_{n+1} = g(a_n)$: the monotone staircase converges to $a^* = 0.8217$ at $d_m = 0.046$ and escapes to $1$ through the gap above the diagonal at $d_m = 0.052$.}
\label{fig:geom}
\end{figure}

\section{Immediate Corollaries}
\label{sec:cor}

\begin{corollary}[Critical surviving fraction jump]
At $d_m = \dc(a_0)$, the fraction of nodes surviving after the cascade completes is strictly positive:
\[
  \prod_{k=0}^{\infty} (1 - p_k) \;=\; 1 - e^{-\frac{1 - a^*}{\dc}} \;=\; 1 - \frac{\dc}{1 + \dc} \;=\; \frac{1}{1 + \dc(a_0)},
\]
whereas for every $d_m > \dc(a_0)$, the surviving fraction drops discontinuously to $0$ (first-order transition).
\end{corollary}
\begin{proof}
Take $n \to \infty$ in \eqref{eq:surv-tel} using $a^* = a_{\mathrm{peak}} = 1 - \dc \ln(1 + 1/\dc)$ from Lemma~\ref{lem:mono}(b), so $e^{-(1-a^*)/\dc} = \frac{\dc}{1+\dc}$.
\end{proof}

\begin{corollary}[Arbitrary uniform capacity $L_{\max} > 0$]
\label{cor:scale}
If every generator has capacity $L_{\max} > 0$ and initial load $a_0 \in [0, L_{\max})$, rescaling loads by $L_{\max}$ leaves the threshold condition dimensionless and yields
\[
  L_{\max} - a_0 \;=\; \dc \left(1 + \ln\!\left(1 + \frac{L_{\max}}{\dc}\right)\right), \qquad
  \dc(a_0; L_{\max}) \;=\; L_{\max}\cdot \dc\!\left(\frac{a_0}{L_{\max}}\right).
\]
\end{corollary}

\begin{corollary}[Exact area under the $\dc(a_0)$ curve]
The area underneath $\dc : [0, 1] \to [0, \dc(0)]$ is
\[
  \int_0^1 \dc(a_0)\,da_0 \;=\; \frac{1}{2}\,\ln\!\bigl(1 + \dc(0)\bigr) \;=\; 0.191248771449\ldots
\]
\end{corollary}
\begin{proof}
Since $\dc$ is strictly decreasing and continuously invertible via $a_0(d) = 1 - d - d\ln(1+1/d)$ from $d = 0$ ($a_0 = 1$) to $d_0 := \dc(0)$ ($a_0 = 0$), swap axes of integration:
\[
  \int_0^1 \dc(a_0)\,da_0 \;=\; \int_0^{d_0} a_0(d)\,dd
  \;=\; \int_0^{d_0} \bigl(1 - d - d\ln(1 + 1/d)\bigr)\,dd.
\]
Using the antiderivative $\frac{d}{dd}\bigl[\tfrac{1}{2}\bigl(d^2\ln(1+1/d) + d - \ln(1+d)\bigr)\bigr] = d\ln(1+1/d)$, we evaluate on $[0, d_0]$ and simplify with $d_0\ln(1+1/d_0) = 1 - d_0$ (i.e.\ $d_0^2\ln(1+1/d_0) = d_0 - d_0^2$):
\[
  \int_0^{d_0} a_0(d)\,dd
  \;=\; d_0 - \frac{d_0^2}{2} - \frac{1}{2}\Bigl[(d_0 - d_0^2) + d_0 - \ln(1 + d_0)\Bigr]
  \;=\; \frac{1}{2}\,\ln(1 + d_0). \qedhere
\]
\end{proof}

\newpage
\section{Numerical Verification and Note on the 2010 Figures}
\label{sec:num}

Table~\ref{tab:steps} directly checks Lemma~\ref{lem:telescope} at $(a_0, d_m) = (0.80, 0.048)$: computing $(a_n, \mu_{n-1}, \md_n, p_n)$ from \eqref{eq:init}--\eqref{eq:rec} agrees with the scalar map $g(a_n)$ to machine precision ($\le 1.1\times 10^{-16}$ across all rounds). Table~\ref{tab:dcrit} lists exact values of $\dc(a_0)$ from \eqref{eq:lambert}.

\vspace{4pt}
\begingroup\centering
\phantomsection
\captionof{table}{Term-by-term comparison of \eqref{eq:init}--\eqref{eq:rec} with the 1D reduction $a_{n+1} = g(a_n)$ at $a_0 = 0.80$, $d_m = 0.048$. Limiting load is $a^* = 0.8309841098$.}
\label{tab:steps}
\small\renewcommand{\arraystretch}{0.91}
\begin{tabular}{@{}rrrrrr@{}}
\toprule
$n$ & $a_n$ & $\mu_{n-1}$ & $\md_n$ & $p_n$ & $|(a_n + \md_n) - g(a_n)|$ \\
\midrule
1 & $0.80000000$ & $1.048000^{\,*}$ & $1.6504\times 10^{-2}$ & $6.4621\times 10^{-3}$ & $0.0\times 10^{0}$ \\
2 & $0.81650391$ & $1.007780$ & $6.5547\times 10^{-3}$ & $3.2709\times 10^{-3}$ & $1.1\times 10^{-16}$ \\
3 & $0.82305860$ & $1.003203$ & $3.2921\times 10^{-3}$ & $1.8252\times 10^{-3}$ & $1.1\times 10^{-16}$ \\
4 & $0.82635073$ & $1.001627$ & $1.8315\times 10^{-3}$ & $1.0729\times 10^{-3}$ & $1.1\times 10^{-16}$ \\
5 & $0.82818224$ & $1.000910$ & $1.0750\times 10^{-3}$ & $6.4976\times 10^{-4}$ & $1.1\times 10^{-16}$ \\
10 & $0.83071515$ & $1.000078$ & $9.8735\times 10^{-5}$ & $6.2371\times 10^{-5}$ & $0.0\times 10^{0}$ \\
15 & $0.83095660$ & $1.000008$ & $1.0055\times 10^{-5}$ & $6.3790\times 10^{-6}$ & $1.1\times 10^{-16}$ \\
25 & $0.83098382$ & $1.000000$ & $1.0655\times 10^{-7}$ & $6.7624\times 10^{-8}$ & $1.1\times 10^{-16}$ \\
40 & $0.83098411$ & $1.000000$ & $1.1643\times 10^{-10}$ & $7.3901\times 10^{-11}$ & $1.1\times 10^{-16}$ \\
\bottomrule
\multicolumn{6}{@{}l}{\scriptsize $^*$At $n=1$, $\mu_0 = 1 + d_m = 1.048$ and $p_0 = e^{-0.2/0.048} = 1.5504\times 10^{-2}$.}
\end{tabular}\par
\endgroup

\vspace{7pt}
\begingroup\centering
\phantomsection
\captionof{table}{Selected exact values of $\dc(a_0)$ from \eqref{eq:lambert}, limiting load $a^* = a_0 + \dc$, surviving fraction $(1+\dc)^{-1}$, and comparison against the sign-flipped $(+)$-bug values appearing in Figs.~5--7 of~\cite{ks2010}.}
\label{tab:dcrit}
\small\renewcommand{\arraystretch}{0.91}
\begin{tabular}{@{}rrrrrr@{}}
\toprule
$a_0$ & Exact $\dc(a_0)$ & $a^* = a_0 + \dc$ & Surviving frac.\ $\tfrac{1}{1+\dc}$ & $(+)$-bug threshold & Residual of \eqref{eq:implicit} \\
\midrule
$0.00$ & $0.465941272385$ & $0.46594127$ & $0.682156$ & $0.393108$ ($-15.6\%$) & $0.0\times 10^{0}$ \\
$0.01$ & $0.459126812951$ & $0.46912681$ & $0.685341$ & $0.388101$ ($-15.5\%$) & $0.0\times 10^{0}$ \\
$0.20$ & $0.336129355299$ & $0.53612936$ & $0.748431$ & $0.295103$ ($-12.2\%$) & $-1.1\times 10^{-16}$ \\
$0.50$ & $0.171002214031$ & $0.67100221$ & $0.853969$ & $0.159722$ ($-6.6\%$) & $-1.1\times 10^{-16}$ \\
\rowcolor{greentint}
$\mathbf{0.80}$ & $\mathbf{0.049281593062}$ & $\mathbf{0.84928159}$ & $\mathbf{0.953033}$ & $0.048460$ ($-1.7\%$) & $0.0\times 10^{0}$ \\
$0.90$ & $0.020345784230$ & $0.92034578$ & $0.980060$ & $0.020227$ ($-0.6\%$) & $-1.4\times 10^{-17}$ \\
$0.95$ & $0.008689097104$ & $0.95868910$ & $0.991386$ & $0.008671$ ($-0.2\%$) & $+6.9\times 10^{-18}$ \\
$0.99$ & $0.001308924995$ & $0.99130892$ & $0.998693$ & $0.001309$ ($-0.02\%$) & $0.0\times 10^{0}$ \\
\bottomrule
\end{tabular}\par
\endgroup

\vspace{10pt}
\begingroup\centering
\phantomsection
\begin{tikzpicture}
  \begin{axis}[
    name=dcplot,
    width=0.45\linewidth, height=4.25cm,
    xmin=0, xmax=1.0, ymin=0, ymax=0.56,
    xlabel={Initial load $a_0$},
    ylabel={$\dc(a_0)$},
    title={\small\sffamily\bfseries (a) $\dc$ vs.\ $a_0$ (cf.\ Fig.~7 of~\cite{ks2010})},
    grid=major, major grid style={draw=rulecol!65,line width=0.3pt},
    axis line style={draw=muted,line width=0.55pt},
    tick label style={font=\scriptsize\sffamily},
    label style={font=\scriptsize\sffamily\color{primary}},
    legend style={font=\tiny\sffamily, at={(0.98,0.97)}, anchor=north east, fill=white, fill opacity=0.92, draw=rulecol, cells={anchor=west}},
    enlargelimits=false
  ]
    \addplot[dotted, thick, color=muted, domain=0.44:1.0] {1 - x};
    \addlegendentry{Excess capacity $1 - a_0$}
    \addplot[thick, line width=1.5pt, color=okgreen] table[x=a0, y=dc] {data/dcrit_curve.dat};
    \addlegendentry{Exact $\dc(a_0)$ \eqref{eq:lambert}}
    \addplot[dashed, thick, line width=1.2pt, color=warnamber] table[x=a0, y=dbug] {data/dcrit_curve.dat};
    \addlegendentry{Flipped $(+)$ in \eqref{eq:rec-mu}}
    \addplot[only marks, mark=*, mark size=1.8pt, color=okgreen] coordinates {(0,0.46594) (0.2,0.33613) (0.5,0.17100) (0.8,0.04928)};
    \addplot[only marks, mark=*, mark size=1.6pt, color=warnamber] coordinates {(0,0.39311) (0.2,0.29510) (0.5,0.15972) (0.8,0.04846)};
    \node[font=\tiny\sffamily, color=okgreen, anchor=south west] at (axis cs:0.02,0.466) {$0.4659$};
    \node[font=\tiny\sffamily, color=warnamber, anchor=north west, fill=white, fill opacity=0.85, text opacity=1, inner sep=0.6pt] at (axis cs:0.015,0.362) {$0.3931$};
    \node[font=\tiny\sffamily, color=okgreen, anchor=south] at (axis cs:0.80,0.068) {$0.0493$};
  \end{axis}

  \begin{axis}[
    at={($(dcplot.east)+(1.2cm,0)$)}, anchor=west,
    width=0.49\linewidth, height=4.25cm,
    xmin=1, xmax=76, ymin=0.795, ymax=1.01,
    xtick={1,10,20,30,40,50,60},
    xlabel={Iteration $n$},
    ylabel={Shifted load $a_n$},
    title={\small\sffamily\bfseries (b) $\{a_n\}$ vs.\ $n$ at $a_0 = 0.80$ (cf.\ Fig.~5)},
    grid=major, major grid style={draw=rulecol!65,line width=0.3pt},
    axis line style={draw=muted,line width=0.55pt},
    tick label style={font=\scriptsize\sffamily},
    label style={font=\scriptsize\sffamily\color{primary}},
    enlargelimits=false
  ]
    \draw[color=muted, densely dotted, line width=0.7pt] (axis cs:1,0.84928) -- (axis cs:60,0.84928);
    \addplot[thick, color=slate!70!white] table[x=n, y=a] {data/an_040.dat};
    \addplot[thick, color=slate] table[x=n, y=a] {data/an_046.dat};
    \addplot[thick, color=okgreen!75!black] table[x=n, y=a] {data/an_048.dat};
    \addplot[thick, line width=1.2pt, color=okgreen] table[x=n, y=a] {data/an_049.dat};
    \addplot[thick, line width=1.2pt, color=crimson] table[x=n, y=a] {data/an_0494.dat};
    \addplot[thick, color=crimson!80!black] table[x=n, y=a] {data/an_050.dat};
    \addplot[thick, color=crimson!60!black] table[x=n, y=a] {data/an_052.dat};
    \addplot[dashed, thick, color=warnamber] table[x=n, y=a] {data/an_bug049.dat};
    \node[font=\fontsize{4.6}{5}\selectfont\sffamily, color=slate!85!black, anchor=west, inner sep=0.4pt] at (axis cs:61,0.807) {$0.040$};
    \node[font=\fontsize{4.6}{5}\selectfont\sffamily, color=slate, anchor=west, inner sep=0.4pt] at (axis cs:61,0.820) {$0.046$};
    \node[font=\fontsize{4.6}{5}\selectfont\sffamily, color=okgreen!75!black, anchor=west, inner sep=0.4pt] at (axis cs:61,0.831) {$0.048$};
    \node[font=\fontsize{4.6}{5}\selectfont\sffamily, color=okgreen, anchor=west, inner sep=0.4pt] at (axis cs:61,0.842) {$0.049$};
    \node[font=\fontsize{4.6}{5}\selectfont\sffamily, color=muted, fill=white, fill opacity=0.95, text opacity=1, anchor=west, inner sep=0.5pt] at (axis cs:61,0.857) {$a^*\!=\!0.849$};
    \node[font=\tiny\sffamily, color=crimson, anchor=east] at (axis cs:46.2,0.960) {$0.0494$};
    \node[font=\tiny\sffamily, color=crimson, anchor=east] at (axis cs:17.2,0.948) {$0.050$};
    \node[font=\tiny\sffamily, color=crimson!75!black, anchor=east] at (axis cs:8.3,0.905) {$0.052$};
    \node[font=\tiny\sffamily, color=warnamber, anchor=west, align=left] at (axis cs:22.8,0.976) {$0.049$ $(+)$-bug\\[-1pt]{\fontsize{4.8}{5}\selectfont hits $1$ at $n\!=\!22$}};
  \end{axis}
\end{tikzpicture}

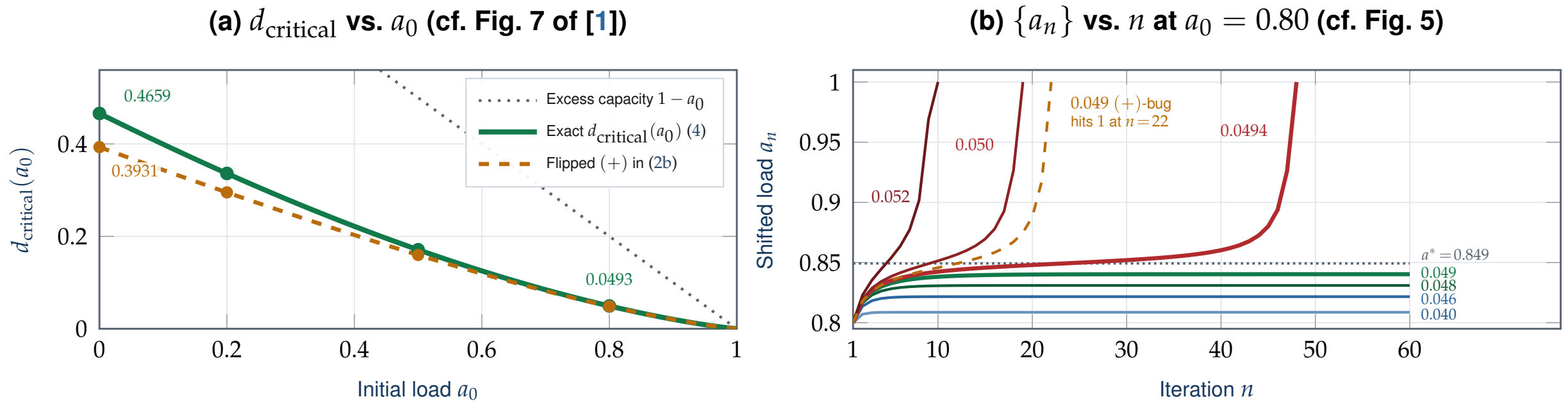
\captionof{figure}{(a) Exact critical threshold $\dc(a_0)$ from \eqref{eq:lambert} compared to excess capacity $1 - a_0$ and the curve produced by flipping $-$ to $+$ in \eqref{eq:rec-mu}. (b) Sample $\{a_n\}$ trajectories from the exact recursion \eqref{eq:init}--\eqref{eq:rec} at $a_0 = 0.80$. Note that $d_m = 0.049$ settles to $0.8403$ under the exact recursion, whereas the sign-flipped $(+)$ version (dashed amber) passes $a_{21} = 0.9175$ and crosses $1$ at $n = 22$.}
\label{fig:curves}\par
\endgroup

\vspace{7pt}
\begin{noteBox}[Remark on the plotted threshold in Figs.~5--7 of \texttt{arXiv:1011.4098v1}]
In \texttt{0802.tex}, equation~\eqref{eq:rec-mu} is printed with the \emph{correct} minus sign, $\mu_{n-1} = 1 + d_m - \frac{\md_{n-1}}{e^{\md_{n-1}/d_m}-1}$. However, the embedded figures in \cite{ks2010} report $\dc \approx 0.048$ at $a_0 = 0.80$ (specifically $0.048460$) rather than $0.049282$. Replacing $-$ by $+$ inside $\mu_{n-1}$ explains the plotted curves quantitatively across all starting loads (Table~\ref{tab:dcrit} and Figure~\ref{fig:curves}): $d_m = 0.049$ stops at $a_{21} = 0.9175$ (matching the last marker before $1$ in Fig.~5 of~\cite{ks2010}) under the flipped sign, while surviving indefinitely at $a^* = 0.840313$ under the exact recursion. Direct $N = 100{,}000$ finite-network Monte Carlo simulations confirm the exact threshold ($a_0 = 0.80$: $d_m = 0.0488$ survives $3/3$ with $\approx 96.4\%$ alive, $d_m = 0.0497$ fails $3/3$; $a_0 = 0.20$: $d_m = 0.310 > 0.295$ survives $3/3$ with $\approx 86.8\%$ alive).
\end{noteBox}

\end{document}